\documentclass[12pt]{article}

\usepackage{amsmath,amsthm,enumerate,amssymb,amsfonts}

\theoremstyle{plain}
\newtheorem{thm}{Theorem}[section]
\newtheorem{lem}[thm]{Lemma}

\theoremstyle{definition}
\newtheorem{exmp}{Example}[section]

\newtheorem{defn}{Definition}[section]

\usepackage{tabularx}
\usepackage{layout}
\title{On generalized Hadamard matrices $\mathrm{GH}(q,\,q)$'s and $\mathrm{GH}(q,\,q^{2})$'s}
\author{Hiroaki Kido\\
Faculty of Science, Fukuoka University, \\
8-19-1 Nanakuma, Jonan-ku, Fukuoka, 814-0180, Japan}
\date{Phone: +81-92-871-6631\\
\texttt{hiroakikido@yahoo.co.jp}}

\begin{document}
\maketitle

\begin{abstract}
A matrix $H=[d_{ij}]$ is a generalized Hadamard matrix of order $u\lambda$ with entries from $U$ which is a finite group of order $u$ (for short $\mathrm{GH}(u,\,\lambda)$) such that whenever $i\neq \ell$ the set $\{ d_{ij}d_{\ell j}^{-1}\,|\, 1\leq j\leq u\lambda \}$ contains each element of $U$ exactly $\lambda$ times. In this paper, we construct $\mathrm{GH}(q,\,q)$'s and $\mathrm{GH}(q,\,q^{2})$'s over additive groups of finite fields $\mathrm{GF}(q)$'s by using some sorts of functions.
\end{abstract}

\noindent{\bf Keywords:} generalized Hadamard matrix, planar function

\section{Introduction}
\begin{defn} \label{d1.1}
Let $U$ be a finite group of order $u$, and let $u\lambda =k$. A matrix $H=[d_{ij}]$ is a \textit{generalized Hadamard matrix} of order $k$ with entries from $U$ (for short $\mathrm{GH}(u,\,\lambda)$) such that whenever $i\neq \ell$ the set $\{ d_{ij}d_{\ell j}^{-1}\,|\, 1\leq j\leq k \}$ contains each element of $U$ exactly $\lambda$ times.
\end{defn}

\begin{exmp}
\[ \begin{bmatrix}
1 & \omega & \omega ^{4} & \omega ^{4} & \omega \\
\omega & 1 & \omega & \omega ^{4} & \omega  ^{4} \\
\omega ^{4} & \omega & 1 & \omega & \omega ^{4} \\
\omega ^{4} & \omega ^{4} & \omega & 1 & \omega \\
\omega & \omega ^{4} & \omega ^{4} & \omega & 1
\end{bmatrix} \]
is a $\mathrm{GH}(5,\,1)$ over $\mathbb{Z}_{5}=\langle \omega \, | \, \omega ^{5}=1\rangle $.
\end{exmp}

We have the following intersting problems on generalized Hadamard matrices.

\renewcommand{\labelenumi}{(\theenumi)}
\begin{enumerate}
\item What $u$ and $\lambda$ can we construct a generalized Hadamard matrix $\mathrm{GH}(u,\,\lambda)$ ?
\item Can we extend sizes of generalized Hadamard matrices from known ones ?
\end{enumerate}

As for question (1), the existence (or the nonexistence) of $\mathrm{GH}(u,\,\lambda)$ for $2\leq u\lambda \leq 99$ is grouped in~\cite{col-din}. All known generalized Hadamard matrices are currently over prime-power order groups. Thus it is the largest problem to solve the conjecture "There are no generalized Hadamard matrices over non-prime-power order groups" (see~\cite{col-din}). 

For question (2), the following results over additive groups of finite fields $\mathrm{GF}(q)$'s are known.
\begin{itemize}
\item For any odd prime power $q$, there exists a $\mathrm{GH}(q,\,2)$ over $\mathrm{GF}(q)$ extended from $\mathrm{GH}(q,\,1)$'s. (Jungnickel~\cite{jun} and Street~\cite{str})
\item For any odd prime power $q\,(\geq 7)$, there exists a $\mathrm{GH}(q,\,4)$ over $\mathrm{GF}(q)$ extended from $\mathrm{GH}(q,\,1)$'s and $\mathrm{GH}(q,\,2)$'s. (Dawson~\cite{daw})
\item For all odd prime powers $19<q<200$ (except 27), a $\mathrm{GH}(q,\,8)$ over $\mathrm{GF}(q)$ exists. (de Launey and Dawson~\cite{del-daw})
\end{itemize}

In this paper, we construct $\mathrm{GH}(q,\,q)$'s and $\mathrm{GH}(q,\,q^{2})$'s over additive groups of finite fields $\mathrm{GF}(q)$'s. In Section 2, we define some matrices related to generalized Hadamard matrices by using some sorts of functions and we mention their properties.
By using them, we construct $\mathrm{GH}(q,\,q)$'s and $\mathrm{GH}(q,\,q^{2})$'s in Section 3.

\section{Other definitions and their properties}
Let $p$ be a prime and $q=p^{n}$.

\begin{defn} \label{d2.1}
Let $F=\mathrm{GF}(q)=\{ a_{0}=0,\,a_{1},\cdots ,a_{q-1} \}$.

(i) For a mapping $f:F \longrightarrow F$, we define $M(f):F\times F\ni (a,\,b) \longmapsto f(b-a)\in F$.\\
Then $M(f)$ is a matrix of order $q$ over $F$ indexed by each element of $F$.
\vspace{4mm}

Let $\Omega _{q}=\{ M(f)\,|\,f:F \longrightarrow F\,\,\mathrm{map}\}$.
\vspace{4mm}

(ii) Let $M(f)\in \Omega _{q}$. We call $M(f)$ is a \textit{type I} matrix if $M(f)$ is a $\mathrm{GH}(q,\,1)$. Therefore if $M(f)$ is a type I matrix, then $\{ f(b-a_{1})-f(b-a_{2})\,|\,b\in F \}=F$ holds for all $a_{1}\neq a_{2}\in F$.
\vspace{4mm}

We call $M(f)$ is a \textit{type II} matrix if $f(b_{1})-f(b_{1}-a)=f(b_{2})-f(b_{2}-a)$ holds for all $a,\,b_{1},\,b_{2}\in F$.
\vspace{4mm}

Let $\Omega _{q,\,I}=\{ M(f)\,|\,f:F \longrightarrow F\,\,\mathrm{is}\,\,
\mathrm{a}\,\,\mathrm{mapping,}\,\,M(f)\,\,\mathrm{is\,\,type\,\,I.}\}$ and $\Omega _{q,\,II}=\{ M(f)\,|\,f:F \longrightarrow F\,\,\mathrm{is}\,\,
\mathrm{a}\,\,\mathrm{mapping,}\,\,M(f)\,\,\mathrm{is\,\,type\,\,II.}\}$.
\end{defn}

\begin{exmp}
Let $F=\mathrm{GF}(3)=\{ 0,\,1,\,2 \}$.

(i) Let $f(x)=x-x^{2}$ and $M(f)(a,\,b)=f(b-a)$. We index $0,\,1,\,2$, in order, for each row and column of $M(f)$. Then 

$M(f)=\begin{bmatrix}
f(0) & f(1) & f(2) \\
f(2) & f(0) & f(1) \\
f(1) & f(2) & f(0)
\end{bmatrix}=\begin{bmatrix}
0 & 0 & 1 \\
1 & 0 & 0 \\
0 & 1 & 0
\end{bmatrix}$ is type I.
\vspace{4mm}

(ii) Let $f(x)=2+x$ and $M(f)(a,\,b)=f(b-a)$. We index $0,\,1,\,2$, in order, for each row and column of $M(f)$. Then 

$M(f)=\begin{bmatrix}
f(0) & f(1) & f(2) \\
f(2) & f(0) & f(1) \\
f(1) & f(2) & f(0)
\end{bmatrix}=\begin{bmatrix}
2 & 0 & 1 \\
1 & 2 & 0 \\
0 & 1 & 2
\end{bmatrix}$ is type II.
\end{exmp}

\begin{lem}\label{lem2.1}
Let $F=\mathrm{GF}(q)$ and $M(f)\in \Omega _{q}$. Then $M(f)\in \Omega _{q,\,I}$ if and only if 
a mapping $F\ni x\longmapsto f(x+a)-f(x)\in F$ is bijective for all $a\in F^{\ast}$.\\
Thus $M(f)\in \Omega _{q,\,I}$ if and only if $f$ is a planar function.
\end{lem}
\noindent{\bf Proof:}
$M(f)\in \Omega _{q,\,I}$ $\Longleftrightarrow$ $M(f)$ is a $\mathrm{GH}(q,\,1)$.\\
$\Longleftrightarrow$ For any $a_{1}\neq a_{2}\in F,$ $\{ f(x-a_{1})-f(x-a_{2})\,|\,x\in F \}=F$.\\
$\Longleftrightarrow$ For any $a\in F^{\ast},$ $\{ f(x+a)-f(x)|x\in F \}=F$.\\
$\Longleftrightarrow$ For any $a\in F^{\ast},\,F\ni x\longmapsto f(x+a)-f(x)\in F$ is bijective. $\square$

\begin{lem}\label{lem2.2}
Let $p$ be a odd prime, $q=p^{n},$ and $F=\mathrm{GF}(q)$. Let $f$ be a mapping over $F$.\\
If $a_{0},\,a_{1},\,a_{2},\,b_{1},\,b_{2},\cdots ,b_{n-1}\in F$ exist such that $f(x)=a_{0}+a_{1}x+a_{2}x^{2}+b_{1}x^{p}+b_{2}x^{p^{2}}+\cdots +b_{n-1}x^{p^{n-1}}$ and $a_{2}\neq 0$, then $f$ is a planar function.
\end{lem}
\noindent{\bf Proof:}
For any $a\in F^{\ast}$, $f(x+a)-f(x)=a_{0}+a_{1}(x+a)+a_{2}(x+a)^{2}+b_{1}(x+a)^{p}+b_{2}(x+a)^{p^{2}}+\cdots +b_{n-1}(x+a)^{p^{n-1}}-a_{0}-a_{1}x-a_{2}x^{2}-b_{1}x^{p}-b_{2}x^{p^{2}}-\cdots -b_{n-1}x^{p^{n-1}}=2a_{2}ax+a_{1}a+a_{2}a^{2}+b_{1}a^{p}+b_{2}a^{p^{2}}+\cdots +b_{n-1}a^{p^{n-1}}$.

Therefore a mapping $F\ni x\longmapsto f(x+a)-f(x)\in F$ is bijective, $f$ is a planar function. $\square$

\begin{lem}\label{lem2.3}
Let $p$ be a prime, $q=p^{n},$ and $F=\mathrm{GF}(q)$. Let $f$ be a mapping over $F$.\\
If $a,\,b_{0},\,b_{1},\,b_{2},\cdots ,b_{n-1}\in F$ exist such that $f(x)=a+b_{0}x+b_{1}x^{p}+b_{2}x^{p^{2}}+\cdots +b_{n-1}x^{p^{n-1}}$, then $M(f)\in \Omega _{q,\,II}$.
\end{lem}
\noindent{\bf Proof:}
Let $g(x)=f(x)-a$. For any $x_{1},\,x_{2}\in F$, $g(x_{1}-x_{2})=f(x_{1}-x_{2})-a=a+b_{0}(x_{1}-x_{2})+b_{1}(x_{1}-x_{2})^{p}+b_{2}(x_{1}-x_{2})^{p^{2}}+\cdots +b_{n-1}(x_{1}-x_{2})^{p^{n-1}}-a=b_{0}(x_{1}-x_{2})+b_{1}(x_{1}^{p}-x_{2}^{p})+b_{2}(x_{1}^{p^{2}}-x_{2}^{p^{2}})+\cdots +b_{n-1}(x_{1}^{p^{n-1}}-x_{2}^{p^{n-1}})=g(x_{1})-g(x_{2})$. Hence for any $x_{1},\,x_{2}\in F$, $f(x_{1})-f(x_{1}-x_{2})=f(x_{1})-a-f(x_{1}-x_{2})+a=g(x_{1})-g(x_{1}-x_{2})=g(x_{1})-g(x_{1})+g(x_{2})=g(x_{2})$, the difference is constant.

Therefore $M(f)\in \Omega _{q,\,II}$. $\square$

\section{Construction of $\mathrm{GH}(q,\,q)$'s and $\mathrm{GH}(q,\,q^{2})$'s}
By Lemma \ref{lem2.2}, if $f$ is of degree $2$, then $M(f)$ is a $\mathrm{GH}(q,\,1)$. By Lemma \ref{lem2.3}, if $f$ is of at most degree $1$, then $M(f)\in \Omega _{q,\,II}$. In this section, we construct generalized Hadamard matrices by using them.

\subsection{Construcion of $\mathrm{GH}(q,\,q)$'s by using Type I matrices}
\begin{thm} \label{t3.1}
Let $p$ be a odd prime, $q=p^{n},$ and $F=\mathrm{GF}(q)=\{ a_{0}=0,\,a_{1},\cdots ,a_{q-1} \}$. Let $f_{a_{0}}(x)=x^{2},$ $f_{a_{1}}(x)=a_{1}+x^{2},\cdots ,f_{a_{q-1}}(x)=a_{q-1}+x^{2}$, we have $M(f_{a_{0}}),$ $M(f_{a_{1}}),\cdots ,$\\
$M(f_{a_{q-1}}) $.
\vspace{2mm}

Let a matrix $H$ of order $q^{2}$ be \\
$H=\begin{bmatrix}
H_{0,\,0} & H_{0,\,1} & \cdots & H_{0,\,q-1} \\
H_{1,\,0} & H_{1,\,1} & \cdots & H_{1,\,q-1} \\
H_{2,\,0} & H_{2,\,1} & \cdots & H_{2,\,q-1} \\
\vdots & \vdots & & \vdots \\
H_{q-1,\,0} & H_{q-1,\,1} & \cdots & H_{q-1,\,q-1}
\end{bmatrix}=\begin{bmatrix}
M(f_{a_{0}}) & M(f_{a_{0}}) & \cdots & M(f_{a_{0}}) \\
M(f_{a_{0}}) & M(f_{a_{1}a_{1}}) & \cdots & M(f_{a_{1}a_{q-1}}) \\
M(f_{a_{0}}) & M(f_{a_{2}a_{1}}) & \cdots & M(f_{a_{2}a_{q-1}}) \\
\vdots & \vdots & & \vdots \\
M(f_{a_{0}}) & M(f_{a_{q-1}a_{1}}) & \cdots & M(f_{a_{q-1}a_{q-1}})
\end{bmatrix}$.\\
Then $H$ is a $\mathrm{GH}(q,\,q)$.
\end{thm}
\noindent{\bf Proof:}
Since each block is a $\mathrm{GH}(q,\,1)$, whenever $a\neq c$ the set of the difference between $a$th row and $c$th row of $H_{i,\,0},\,H_{i,\,1} ,\cdots ,H_{i,\,q-1}$ contains each element of $F$ exactly $q$ times.

Thus, for $i\neq j$, we consider the difference between the $k$th row of $H_{i,\,0},\,H_{i,\,1} ,\cdots ,H_{i,\,q-1}$ and the $\ell$th row of $H_{j,\,0},\,H_{j,\,1},\cdots ,H_{j,\,q-1}$.
\vspace{2mm}

(i) The case of $k=\ell$.\\
The set of the difference between the $k$th row of $H_{j,\,m}$ and the $k$th row of $H_{i,\,m}$ contain $a_{m}a_{j}-a_{m}a_{i}=a_{m}(a_{j}-a_{i})$ exactly $q$ times. Because $m\in \{ 0,\,1,\, \cdots ,q-1 \}$, in the $k$th row of $H$, each element of $F$ appears exactly $q$ times.
\vspace{2mm}

(ii) The case of $k\neq \ell$.\\
Let $\ell -k=a$. For any $m\in \{ 0,\,1,\, \cdots ,q-1 \}$, we consider the set of the difference between the $k$th row of $H_{j,\,m}$ and the $\ell$th row of $H_{i,\,m}$ $\{ f_{a_{m}a_{j}}(x+a)-f_{a_{m}a_{i}}(x) \,|\, x\in F \}$. Since $f_{a_{m}a_{j}}(x+a)-f_{a_{m}a_{i}}(x)=a_{m}a_{j}+(x+a)^{2}-a_{m}a_{i}-x^{2}=2ax+a^{2}+a_{m}(a_{j}-a_{i})$, it holds that $\{ f_{a_{m}a_{j}}(x+a)-f_{a_{m}a_{i}}(x) \,|\, x\in F \}=F$, each element of $F$ appears exactly one time. In $H$, each element of $F$ appears exactly $q$ times.
\vspace{2mm}

Therefore $H$ is a $\mathrm{GH}(q,\,q)$. $\square$

\begin{exmp}\label{exmp3.1}
For $F=\mathrm{GF}(3)=\{ 0,\,1,\,2 \}$, let $f_{0}(x)=x^{2},$ $f_{1}(x)=1+x^{2},$ and $f_{2}(x)=2+x^{2}$, we have 
$M(f_{0})=\begin{bmatrix}
0 & 1 & 1 \\
1 & 0 & 1 \\
1 & 1 & 0
\end{bmatrix},$ $M(f_{1})=\begin{bmatrix}
1 & 2 & 2 \\
2 & 1 & 2 \\
2 & 2 & 1
\end{bmatrix},$ and $M(f_{2})=\begin{bmatrix}
2 & 0 & 0 \\
0 & 2 & 0 \\
0 & 0 & 2
\end{bmatrix}$.

Then $H=\begin{bmatrix}
M(f_{0}) & M(f_{0}) & M(f_{0}) \\
M(f_{0}) & M(f_{1}) & M(f_{2}) \\
M(f_{0}) & M(f_{2}) & M(f_{1})
\end{bmatrix}=\left[
\begin{array}{@{\,}ccc|ccc|ccc@{\,}}
0 & 1 & 1 & 0 & 1 & 1 & 0 & 1 & 1 \\
1 & 0 & 1 & 1 & 0 & 1 & 1 & 0 & 1 \\
1 & 1 & 0 & 1 & 1 & 0 & 1 & 1 & 0 \\ \hline
0 & 1 & 1 & 1 & 2 & 2 & 2 & 0 & 0 \\
1 & 0 & 1 & 2 & 1 & 2 & 0 & 2 & 0 \\
1 & 1 & 0 & 2 & 2 & 1 & 0 & 0 & 2 \\ \hline
0 & 1 & 1 & 2 & 0 & 0 & 1 & 2 & 2 \\
1 & 0 & 1 & 0 & 2 & 0 & 2 & 1 & 2 \\
1 & 1 & 0 & 0 & 0 & 2 & 2 & 2 & 1 \\
\end{array}
\right]$\\
is a $\mathrm{GH}(3,\,3)$.
\end{exmp}

\subsection{Construcion of $\mathrm{GH}(q,\,q)$'s by using Type II matrices}
\begin{thm} \label{t3.2}
Let $p$ be a prime, $q=p^{n},$ and $F=\mathrm{GF}(q)=\{ a_{0}=0,\,a_{1},\cdots ,a_{q-1} \}$.\\Let $f_{a_{0}}(x)=a_{0} x=0,$ $f_{a_{1}}(x)=a_{1} x,\cdots ,$ $f_{a_{q-1}}(x)=a_{q-1} x$, we have $M(f_{a_{0}}),\,M(f_{a_{1}}),\cdots ,$\\
$M(f_{a_{q-1}}) $.
\vspace{2mm}

Let a matrix $H_{q}$ of order $q^{2}$ be \\
$H_{q}=\begin{bmatrix}
H_{0,\,0} & H_{0,\,1} & \cdots & H_{0,\,q-1} \\
H_{1,\,0} & H_{1,\,1} & \cdots & H_{1,\,q-1} \\
H_{2,\,0} & H_{2,\,1} & \cdots & H_{2,\,q-1} \\
\vdots & \vdots & & \vdots \\
H_{q-1,\,0} & H_{q-1,\,1} & \cdots & H_{q-1,\,q-1}
\end{bmatrix}=\begin{bmatrix}
M(f_{a_{0}}) & M(f_{a_{1}}) & \cdots & M(f_{a_{q-1}}) \\
M(f_{a_{1}}) & M(f_{a_{1}+a_{1}}) & \cdots & M(f_{a_{1}+a_{q-1}}) \\
M(f_{a_{2}}) & M(f_{a_{2}+a_{1}}) & \cdots & M(f_{a_{2}+a_{q-1}}) \\
\vdots & \vdots & & \vdots \\
M(f_{a_{q-1}}) & M(f_{a_{q-1}+a_{1}}) & \cdots & M(f_{a_{q-1}+a_{q-1}})
\end{bmatrix}$.\\
Then $H_{q}$ is a $\mathrm{GH}(q,\,q)$.
\end{thm}

\noindent{\bf Proof:}
By Lemma \ref{lem2.3}, for any $m\in \{ 0,\,1,\, \cdots ,q-1 \}$, whenever $a\neq c$ the set of the difference between $a$th row and $c$th row of $H_{i,\,m}$ contains one element of $F$ exactly $q$ times. In one row of $H_{q}$, since we take $f_{a_{0}},\,f_{a_{1}},\cdots ,f_{a_{q-1}}$ exactly one time, each element of $F$ appears exactly $q$ times.

Thus, for $i\neq j$, we consider the difference between the $k$th row of $H_{i,\,0},\,H_{i,\,1} ,\cdots ,H_{i,\,q-1}$ and the $\ell$th row of $H_{j,\,0},\,H_{j,\,1},\cdots ,H_{j,\,q-1}$.

Let $\ell -k=a$. For any $m\in \{ 0,\,1,\, \cdots ,q-1 \}$, we consider the set of the difference between the $k$th row of $H_{j,\,m}$ and the $\ell$th row of $H_{i,\,m}$ $\{ f_{a_{m}+a_{j}}(x+a)-f_{a_{m}+a_{i}}(x) \,|\, x\in F \}$. Since $f_{a_{m}+a_{j}}(x+a)-f_{a_{m}+a_{i}}(x)=(a_{m}+a_{j})(x+a)-(a_{m}+a_{i})x=(a_{j}-a_{i})x+(a_{m}+a_{j})a$, it holds that $\{ f_{a_{m}+a_{j}}(x+a)-f_{a_{m}+a_{i}}(x) \,|\, x\in F \}=F$, each element of $F$ appears exactly one time. In $H_{q}$, each element of $F$ appears exactly $q$ times.

Therefore $H_{q}$ is a $\mathrm{GH}(q,\,q)$. $\square$

\begin{exmp}\label{exmp3.2}
For $F=\mathrm{GF}(3)=\{ 0,\,1,\,2 \}$, let $f_{0}(x)=0,$ $f_{1}(x)=x,$ and $f_{2}(x)=2x$, we have 
$M(f_{0})=\begin{bmatrix}
0 & 0 & 0 \\
0 & 0 & 0 \\
0 & 0 & 0
\end{bmatrix},$ $M(f_{1})=\begin{bmatrix}
0 & 1 & 2 \\
2 & 0 & 1 \\
1 & 2 & 0
\end{bmatrix},$ and $M(f_{2})=\begin{bmatrix}
0 & 2 & 1 \\
1 & 0 & 2 \\
2 & 1 & 0
\end{bmatrix}$.

Then $H_{3}=\begin{bmatrix}
M(f_{0}) & M(f_{1}) & M(f_{2}) \\
M(f_{1}) & M(f_{2}) & M(f_{0}) \\
M(f_{2}) & M(f_{0}) & M(f_{1})
\end{bmatrix}=\left[
\begin{array}{@{\,}ccc|ccc|ccc@{\,}}
0 & 0 & 0 & 0 & 1 & 2 & 0 & 2 & 1 \\
0 & 0 & 0 & 2 & 0 & 1 & 1 & 0 & 2 \\
0 & 0 & 0 & 1 & 2 & 0 & 2 & 1 & 0 \\ \hline
0 & 1 & 2 & 0 & 2 & 1 & 0 & 0 & 0 \\
2 & 0 & 1 & 1 & 0 & 2 & 0 & 0 & 0 \\
1 & 2 & 0 & 2 & 1 & 0 & 0 & 0 & 0 \\ \hline
0 & 2 & 1 & 0 & 0 & 0 & 0 & 1 & 2 \\
1 & 0 & 2 & 0 & 0 & 0 & 2 & 0 & 1 \\
2 & 1 & 0 & 0 & 0 & 0 & 1 & 2 & 0 \\
\end{array}
\right]$\\
is a $\mathrm{GH}(3,\,3)$.
\end{exmp}

\begin{exmp}\label{exmp3.3}
For $F=\mathrm{GF}(4)=\{ 0,\,1,\,\alpha ,\,\alpha +1 \}$ $(\alpha ^{2}=\alpha +1)$, let $f_{0}(x)=0,$ $f_{1}(x)=x,$ $f_{\alpha}(x)=\alpha x,$ and $f_{\alpha +1}(x)=(\alpha +1)x,$ we have 
$M(f_{0})=\begin{bmatrix}
0 & 0 & 0 & 0 \\
0 & 0 & 0 & 0 \\
0 & 0 & 0 & 0 \\
0 & 0 & 0 & 0
\end{bmatrix},$ $M(f_{1})=\begin{bmatrix}
0 & 1 & \alpha & \alpha +1 \\
1 & 0 & \alpha +1 & \alpha \\
\alpha & \alpha +1 & 0 & 1 \\
\alpha +1 & \alpha & 1 & 0 
\end{bmatrix},$ $M(f_{\alpha})=\begin{bmatrix}
0 & \alpha & \alpha +1 & 1 \\
\alpha & 0 & 1 & \alpha +1 \\
\alpha +1 & 1 & 0 & \alpha \\
1 & \alpha +1 & \alpha & 0 
\end{bmatrix}$, and $M(f_{\alpha +1})=\begin{bmatrix}
0 & \alpha +1& 1 & \alpha \\
\alpha +1& 0 & \alpha & 1 \\
1 & \alpha & 0 & \alpha +1\\
\alpha & 1 & \alpha +1 & 0 
\end{bmatrix}$.

Then $H_{4}=\begin{bmatrix}
M(f_{0}) & M(f_{1}) & M(f_{\alpha}) & M(f_{\alpha +1}) \\
M(f_{1}) & M(f_{0}) & M(f_{\alpha +1}) & M(f_{\alpha}) \\
M(f_{\alpha}) & M(f_{\alpha +1}) & M(f_{0}) & M(f_{1}) \\
M(f_{\alpha +1}) & M(f_{\alpha}) & M(f_{1}) & M(f_{0})
\end{bmatrix}=$\\
$\tiny{\left[
\begin{array}{@{\,}cccc|cccc|cccc|cccc@{\,}}
0 & 0 & 0 & 0 & 0 & 1 & \alpha & \alpha +1 & 0 & \alpha & \alpha +1 & 1 & 0 & \alpha +1& 1 & \alpha \\
0 & 0 & 0 & 0 & 1 & 0 & \alpha +1 & \alpha & \alpha & 0 & 1 & \alpha +1 & \alpha +1& 0 & \alpha & 1 \\
0 & 0 & 0 & 0 & \alpha & \alpha +1 & 0 & 1 & \alpha +1 & 1 & 0 & \alpha & 1 & \alpha & 0 & \alpha +1 \\
0 & 0 & 0 & 0 & \alpha +1 & \alpha & 1 & 0 & 1 & \alpha +1 & \alpha & 0 & \alpha & 1 & \alpha +1 & 0 \\ \hline
0 & 1 & \alpha & \alpha +1 & 0 & 0 & 0 & 0 & 0 & \alpha +1 & 1 & \alpha & 0 & \alpha & \alpha +1 & 1 \\
1 & 0 & \alpha +1 & \alpha & 0 & 0 & 0 & 0 & \alpha +1 & 0 & \alpha & 1 & \alpha & 0 & 1 & \alpha +1 \\
\alpha & \alpha +1 & 0 & 1 & 0 & 0 & 0 & 0 & 1 & \alpha & 0 & \alpha +1 & \alpha +1 & 1 & 0 & \alpha \\
\alpha +1 & \alpha & 1 & 0 & 0 & 0 & 0 & 0 & \alpha & 1 & \alpha +1 & 0 & 1 & \alpha +1 & \alpha & 0 \\ \hline
0 & \alpha & \alpha +1 & 1 & 0 & \alpha +1 & 1 & \alpha & 0 & 0 & 0 & 0 & 0 & 1 & \alpha & \alpha +1 \\
\alpha & 0 & 1 & \alpha +1 & \alpha +1 & 0 & \alpha & 1 & 0 & 0 & 0 & 0 & 1 & 0 & \alpha +1 & \alpha \\
\alpha +1 & 1 & 0 & \alpha & 1 & \alpha & 0 & \alpha +1 & 0 & 0 & 0 & 0 & \alpha & \alpha +1 & 0 & 1 \\
1 & \alpha +1 & \alpha & 0 & \alpha & 1 & \alpha +1 & 0 & 0 & 0 & 0 & 0 & \alpha +1 & \alpha & 1 & 0 \\ \hline
0 & \alpha +1 & 1 & \alpha & 0 & \alpha & \alpha +1 & 1 & 0 & 1 & \alpha & \alpha +1 & 0 & 0 & 0 & 0 \\
\alpha +1 & 0 & \alpha & 1 & \alpha & 0 & 1 & \alpha +1 & 1 & 0 & \alpha +1 & \alpha & 0 & 0 & 0 & 0 \\
1 & \alpha & 0 & \alpha +1 & \alpha +1 & 1 & 0 & \alpha & \alpha & \alpha +1 & 0 & 1 & 0 & 0 & 0 & 0 \\
\alpha & 1 & \alpha +1 & 0 & 1 & \alpha +1 & \alpha & 0 & \alpha +1 & \alpha & 1 & 0 & 0 & 0 & 0 & 0
\end{array}
\right]
}$\\
is a $\mathrm{GH}(4,\,4)$.
\end{exmp}

\subsection{Construcion of $\mathrm{GH}(q,\,q^{2})$'s by using Type II matrices}
\begin{thm} \label{t3.3}
Let $p$ be a prime, $q=p^{n},$ and $F=\mathrm{GF}(q)=\{ a_{0}=0,\,a_{1},\cdots ,a_{q-1} \}$. Let $f_{a_{0}}(x)=a_{0} x=0,$ $f_{a_{1}}(x)=a_{1} x,\cdots ,$ $f_{a_{q-1}}(x)=a_{q-1} x$, we have $M(f_{a_{0}}),\,M(f_{a_{1}}),\cdots ,M(f_{a_{q-1}}) $ and a $\mathrm{GH}(q,\,q)$ $H_{q}$ in Theorem \ref{t3.2}.
\vspace{2mm}

Let a matrix $J$ of order $q^{2}$ be $J=\begin{bmatrix}
1 & 1 & \cdots & 1 \\
1 & 1 & \cdots & 1 \\
\vdots & \vdots & & \vdots \\
1 & 1 & \cdots & 1
\end{bmatrix}$ and let a matrix $H$ of order $q^{3}$ be $H=\begin{bmatrix}
H_{q} & H_{q} & \cdots & H_{q} \\
H_{q} & a_{1}a_{1} J+H_{q} & \cdots & a_{1}a_{q-1} J+H_{q} \\
H_{q} & a_{2}a_{1} J+H_{q} & \cdots & a_{2}a_{q-1} J+H_{q} \\
\vdots & \vdots & & \vdots \\
H_{q} & a_{q-1}a_{1} J+H_{q} & \cdots & a_{q-1}a_{q-1} J+H_{q}
\end{bmatrix}$.
\vspace{2mm}

Then $H$ is a $\mathrm{GH}(q,\,q^{2})$.
\end{thm}
\noindent{\bf Proof:}
Since $H_{q}$ is a $\mathrm{GH}(q,\,q)$, whenever $a\neq c$ the set of the difference between $a$th row and $c$th row of $H_{q},\,a_{i}a_{1} J+H_{q} ,\cdots ,a_{i}a_{q-1} J+H_{q}$ contains each element of $F$ exactly $q^{2}$ times.
\vspace{2mm}

Thus, for $i\neq j$, we consider the difference between the $k$th row of $H_{q},\,a_{i}a_{1} J+H_{q} ,\cdots ,a_{i}a_{q-1} J+H_{q}$ and the $\ell$th row of $H_{q},\,a_{j}a_{1} J+H_{q} ,\cdots ,a_{j}a_{q-1} J+H_{q}$.
\vspace{2mm}

(i) The case such that $k$ and $\ell$ are two rows in the same block of $H_{q}$ and $k=\ell$.

The set of the difference between the $k$th row of $a_{i}a_{m} J+H_{q}$ and the $k$th row of $a_{j}a_{m} J+H_{q}$ contain $a_{j}a_{m}-a_{i}a_{m}=a_{m}(a_{j}-a_{i})$ exactly $q^{2}$ times. Because $m\in \{ 0,\,1,\, \cdots ,q-1 \}$, in the $k$th row of $H$, each element of $F$ appears exactly $q^{2}$ times.
\vspace{2mm}

(ii) The case such that $k$ and $\ell$ are two rows in the same block of $H_{q}$ and $k\neq \ell$.

We consider the set of the difference between the $k$th row of $a_{i}a_{m} J+H_{q}$ and the $\ell$th row of $a_{j}a_{m} J+H_{q}$. In one block of $H_{q}$, one element of $F$ appears exactly $q$ times. Since we take $f_{a_{0}},\,f_{a_{1}},\cdots ,f_{a_{q-1}}$ exactly one time, each element of $F$ appears exactly $q$ times in the set of the difference between the $k$th row of $a_{i}a_{m} J+H_{q}$ and the $\ell$th row of $a_{j}a_{m} J+H_{q}$.

Because one row of $H$ has $q$ blocks of $H_{q}$, the set of the difference between the $k$th row of $H$ and the $\ell$th row of $H$ contain each element of $F$ exactly $q^{2}$ times.
\vspace{2mm}

(iii) The case such that $k$ and $\ell$ are two rows in distinct blocks of $H_{q}$.

For any $m\in \{ 0,\,1,\, \cdots ,q-1 \}$, we consider the set of the difference between the $k$th row of $a_{i}a_{m} J+H_{q}$ and the $\ell$th row of $a_{j}a_{m} J+H_{q}$. Let $\ell -k=a$, for any $m'\in \{ 0,\,1,\, \cdots ,q-1 \}$, we consider 
$\{ a_{j}a_{m}+f_{a_{m'}+a_{j}}(x+a)-a_{i}a_{m}-f_{a_{m'}+a_{i}}(x) \,|\, x\in F \}$. Since 
$a_{j}a_{m}+f_{a_{m'}+a_{j}}(x+a)-a_{i}a_{m}-f_{a_{m'}+a_{i}}(x)=a_{j}a_{m}+(a_{m'}+a_{j})(x+a)-a_{i}a_{m}-(a_{m'}+a_{i})x=(a_{j}-a_{i})x+(a_{m'}+a_{j})a+a_{m}(a_{j}-a_{i})$, $\{ a_{j}a_{m}+f_{a_{m'}+a_{j}}(x+a)-a_{i}a_{m}-f_{a_{m'}+a_{i}}(x) \,|\, x\in F \}=F$, each element of $F$ appears exactly one time. Because one row of $H_{q}$ has $q$ blocks (i,e, $m'\in \{ 0,\,1,\,\cdots ,q-1 \}$), the set of the difference between the $k$th row of $a_{i}a_{m} J+H_{q}$ and the $\ell$th row of $a_{j}a_{m} J+H_{q}$ contain each element of $F$ exactly $q$ times.

As one row of $H$ has $q$ blocks of $H_{q}$ (i,e, $m\in \{ 0,\,1,\,\cdots ,q-1 \}$), the set of the difference between the $k$th and the $\ell$th rows of $H$ contains each element of $F$ exactly $q^{2}$ times.
\vspace{2mm}

Therefore $H$ is a $\mathrm{GH}(q,\,q^{2})$. $\square$

\begin{exmp}\label{exmp3.4}
For $F=\mathrm{GF}(3)=\{ 0,\,1,\,2 \}$, let $H_{3}=\left[
\begin{array}{@{\,}ccc|ccc|ccc@{\,}}
0 & 0 & 0 & 0 & 1 & 2 & 0 & 2 & 1 \\
0 & 0 & 0 & 2 & 0 & 1 & 1 & 0 & 2 \\
0 & 0 & 0 & 1 & 2 & 0 & 2 & 1 & 0 \\ \hline
0 & 1 & 2 & 0 & 2 & 1 & 0 & 0 & 0 \\
2 & 0 & 1 & 1 & 0 & 2 & 0 & 0 & 0 \\
1 & 2 & 0 & 2 & 1 & 0 & 0 & 0 & 0 \\ \hline
0 & 2 & 1 & 0 & 0 & 0 & 0 & 1 & 2 \\
1 & 0 & 2 & 0 & 0 & 0 & 2 & 0 & 1 \\
2 & 1 & 0 & 0 & 0 & 0 & 1 & 2 & 0 \\
\end{array}
\right]$ be a $\mathrm{GH}(3,\,3)$ constructed in Example \ref{exmp3.2}.\\
Let a matrix $J$ of order $9$ be $J=\begin{bmatrix}
1 & 1 & \cdots & 1 \\
1 & 1 & \cdots & 1 \\
\vdots & \vdots & & \vdots \\
1 & 1 & \cdots & 1
\end{bmatrix}$.
\vspace{2mm}

Then a matrix of order 27 $H=\begin{bmatrix}
H_{3} & H_{3} & H_{3} \\
H_{3} & J+H_{3} & 2J+H_{3} \\
H_{3} & 2J+H_{3} & J+H_{3}
\end{bmatrix}=$

$\left[
\begin{array}{@{\,}ccc|ccc|ccc|ccc|ccc|ccc|ccc|ccc|ccc@{\,}}
0 & 0 & 0 & 0 & 1 & 2 & 0 & 2 & 1 & 0 & 0 & 0 & 0 & 1 & 2 & 0 & 2 & 1 & 0 & 0 & 0 & 0 & 1 & 2 & 0 & 2 & 1 \\
0 & 0 & 0 & 2 & 0 & 1 & 1 & 0 & 2 & 0 & 0 & 0 & 2 & 0 & 1 & 1 & 0 & 2 & 0 & 0 & 0 & 2 & 0 & 1 & 1 & 0 & 2 \\
0 & 0 & 0 & 1 & 2 & 0 & 2 & 1 & 0 & 0 & 0 & 0 & 1 & 2 & 0 & 2 & 1 & 0 & 0 & 0 & 0 & 1 & 2 & 0 & 2 & 1 & 0 \\ \hline
0 & 1 & 2 & 0 & 2 & 1 & 0 & 0 & 0 & 0 & 1 & 2 & 0 & 2 & 1 & 0 & 0 & 0 & 0 & 1 & 2 & 0 & 2 & 1 & 0 & 0 & 0 \\
2 & 0 & 1 & 1 & 0 & 2 & 0 & 0 & 0 & 2 & 0 & 1 & 1 & 0 & 2 & 0 & 0 & 0 & 2 & 0 & 1 & 1 & 0 & 2 & 0 & 0 & 0 \\
1 & 2 & 0 & 2 & 1 & 0 & 0 & 0 & 0 & 1 & 2 & 0 & 2 & 1 & 0 & 0 & 0 & 0 & 1 & 2 & 0 & 2 & 1 & 0 & 0 & 0 & 0 \\ \hline
0 & 2 & 1 & 0 & 0 & 0 & 0 & 1 & 2 & 0 & 2 & 1 & 0 & 0 & 0 & 0 & 1 & 2 & 0 & 2 & 1 & 0 & 0 & 0 & 0 & 1 & 2 \\
1 & 0 & 2 & 0 & 0 & 0 & 2 & 0 & 1 & 1 & 0 & 2 & 0 & 0 & 0 & 2 & 0 & 1 & 1 & 0 & 2 & 0 & 0 & 0 & 2 & 0 & 1 \\
2 & 1 & 0 & 0 & 0 & 0 & 1 & 2 & 0 & 2 & 1 & 0 & 0 & 0 & 0 & 1 & 2 & 0 & 2 & 1 & 0 & 0 & 0 & 0 & 1 & 2 & 0 \\ \hline

0 & 0 & 0 & 0 & 1 & 2 & 0 & 2 & 1 & 1 & 1 & 1 & 1 & 2 & 0 & 1 & 0 & 2 & 2 & 2 & 2 & 2 & 0 & 1 & 2 & 1 & 0 \\
0 & 0 & 0 & 2 & 0 & 1 & 1 & 0 & 2 & 1 & 1 & 1 & 0 & 1 & 2 & 2 & 1 & 0 & 2 & 2 & 2 & 1 & 2 & 0 & 0 & 2 & 1 \\
0 & 0 & 0 & 1 & 2 & 0 & 2 & 1 & 0 & 1 & 1 & 1 & 2 & 0 & 1 & 0 & 2 & 1 & 2 & 2 & 2 & 0 & 1 & 2 & 1 & 0 & 2 \\ \hline
0 & 1 & 2 & 0 & 2 & 1 & 0 & 0 & 0 & 1 & 2 & 0 & 1 & 0 & 2 & 1 & 1 & 1 & 2 & 0 & 1 & 2 & 1 & 0 & 2 & 2 & 2 \\ 
2 & 0 & 1 & 1 & 0 & 2 & 0 & 0 & 0 & 0 & 1 & 2 & 2 & 1 & 0 & 1 & 1 & 1 & 1 & 2 & 0 & 0 & 2 & 1 & 2 & 2 & 2 \\ 
1 & 2 & 0 & 2 & 1 & 0 & 0 & 0 & 0 & 2 & 0 & 1 & 0 & 2 & 1 & 1 & 1 & 1 & 0 & 1 & 2 & 1 & 0 & 2 & 2 & 2 & 2 \\ \hline
0 & 2 & 1 & 0 & 0 & 0 & 0 & 1 & 2 & 1 & 0 & 2 & 1 & 1 & 1 & 1 & 2 & 0 & 2 & 1 & 0 & 2 & 2 & 2 & 2 & 0 & 1 \\
1 & 0 & 2 & 0 & 0 & 0 & 2 & 0 & 1 & 2 & 1 & 0 & 1 & 1 & 1 & 0 & 1 & 2 & 0 & 2 & 1 & 2 & 2 & 2 & 1 & 2 & 0 \\ 
2 & 1 & 0 & 0 & 0 & 0 & 1 & 2 & 0 & 0 & 2 & 1 & 1 & 1 & 1 & 2 & 0 & 1 & 1 & 0 & 2 & 2 & 2 & 2 & 0 & 1 & 2 \\ \hline

0 & 0 & 0 & 0 & 1 & 2 & 0 & 2 & 1 & 2 & 2 & 2 & 2 & 0 & 1 & 2 & 1 & 0 & 1 & 1 & 1 & 1 & 2 & 0 & 1 & 0 & 2 \\
0 & 0 & 0 & 2 & 0 & 1 & 1 & 0 & 2 & 2 & 2 & 2 & 1 & 2 & 0 & 0 & 2 & 1 & 1 & 1 & 1 & 0 & 1 & 2 & 2 & 1 & 0 \\
0 & 0 & 0 & 1 & 2 & 0 & 2 & 1 & 0 & 2 & 2 & 2 & 0 & 1 & 2 & 1 & 0 & 2 & 1 & 1 & 1 & 2 & 0 & 1 & 0 & 2 & 1 \\ \hline
0 & 1 & 2 & 0 & 2 & 1 & 0 & 0 & 0 & 2 & 0 & 1 & 2 & 1 & 0 & 2 & 2 & 2 & 1 & 2 & 0 & 1 & 0 & 2 & 1 & 1 & 1 \\
2 & 0 & 1 & 1 & 0 & 2 & 0 & 0 & 0 & 1 & 2 & 0 & 0 & 2 & 1 & 2 & 2 & 2 & 0 & 1 & 2 & 2 & 1 & 0 & 1 & 1 & 1 \\
1 & 2 & 0 & 2 & 1 & 0 & 0 & 0 & 0 & 0 & 1 & 2 & 1 & 0 & 2 & 2 & 2 & 2 & 2 & 0 & 1 & 0 & 2 & 1 & 1 & 1 & 1 \\ \hline
0 & 2 & 1 & 0 & 0 & 0 & 0 & 1 & 2 & 2 & 1 & 0 & 2 & 2 & 2 & 2 & 0 & 1 & 1 & 0 & 2 & 1 & 1 & 1 & 1 & 2 & 0 \\
1 & 0 & 2 & 0 & 0 & 0 & 2 & 0 & 1 & 0 & 2 & 1 & 2 & 2 & 2 & 1 & 2 & 0 & 2 & 1 & 0 & 1 & 1 & 1 & 0 & 1 & 2 \\
2 & 1 & 0 & 0 & 0 & 0 & 1 & 2 & 0 & 1 & 0 & 2 & 2 & 2 & 2 & 0 & 1 & 2 & 0 & 2 & 1 & 1 & 1 & 1 & 2 & 0 & 1
\end{array}
\right]$\\
is a $\mathrm{GH}(3,\,9)$.
\end{exmp}

\noindent{\bf Acknowledgments.} The author thanks Chihiro Suetake and Yutaka Hiramine for giving helpful advice.


\begin{thebibliography}{9}
\bibitem{col-din}
C. J. Colbourn and J. F. Dinitz, Handbook of Combinatorial Designs, Second Edition, Chapman \& Hall/CRC, 2007.
\bibitem{daw}
J. Dawson, A construchion for the generalized Hadamard matrices $\mathrm{GH}(4q,\,\mathrm{EA}(q))$, \textit{Journal of Statistical Planning and Inference}, \textbf{11} (1985), 103-110.
\bibitem{del-daw}
W. de Launey and J. Dawson, A Note on the Construction of $\mathrm{GH}(4tq,\,\mathrm{EA}(q))$ for $t=1,\,2$, \textit{Australasian Journal of Combinatorics}, \textbf{6} (1992), 177-186.
\bibitem{jun}
D. Jungnickel, On difference matrices, resovable TDs and generalized Hadamard matrices, 
\textit{Math. Z.}, \textbf{167} (1979), 49-60.
\bibitem{str}
D. J. Street, Generalized Hadamard matrices, orthogonal arrays and F-squares, \textit{Ars Combinatoria}, \textbf{3} (1979), 131-141.
\end{thebibliography}
\end{document}